\documentclass[11pt,a4paper]{article}
\usepackage[utf8]{inputenc}
\usepackage[T1]{fontenc}
\usepackage[numbers,square]{natbib}
\usepackage{amsmath, amssymb, amsthm}
\usepackage{mathrsfs}
\usepackage{geometry}
\usepackage{hyperref}
\usepackage{enumitem}
\usepackage{booktabs}
\usepackage{graphicx}
\usepackage{xcolor}
\usepackage{bm}
\usepackage[capitalize]{cleveref}
\usepackage{mathtools}
\usepackage{lmodern}

\hypersetup{
    colorlinks=true,
    linkcolor=blue,
    citecolor=red,
    urlcolor=cyan,
    pdftitle={Stochastic Burgers Equation Driven by an Additive Hermite Sheet},
    pdfauthor={Atef Lechiheb}
}

\newtheorem{theorem}{Theorem}[section]
\newtheorem{lemma}[theorem]{Lemma}
\newtheorem{proposition}[theorem]{Proposition}

\newtheorem{definition}[theorem]{Definition}
\newtheorem{remark}[theorem]{Remark}
\newtheorem{example}[theorem]{Example}
\newtheorem{assumption}[theorem]{Assumption}

\newcommand{\E}{\mathbb{E}}

\title{\textbf{Stochastic Burgers Equation Driven by a Hermite Sheet with Additive Noise: \\ Existence, Uniqueness, and Regularity}}
\author{Atef Lechiheb \\
Toulouse School of Economics,\\
Universit\'e Toulouse Capitole, Toulouse, France\\
Email: \href{mailto:atef.lechiheb@gmail.com}{atef.lechiheb@gmail.com}}

\begin{document}

\maketitle

\begin{abstract}

We study the stochastic Burgers equation driven by an additive Hermite sheet of order $q \ge 1$. The equation is formulated in the mild sense using the heat semigroup, and existence and uniqueness of solutions are established via a fixed-point argument in suitable Banach spaces. Under appropriate conditions on the Hurst parameters of the Hermite sheet, we derive uniform moment estimates for the solution, which form the basis for the regularity analysis.

We prove that the solution admits a continuous modification that is H\"older continuous in both time and space, with exponents determined by the Hurst parameters of the driving noise. In addition, we show that the solution inherits an anisotropic self-similarity property from the Hermite sheet, and we identify the corresponding scaling exponents.

The additive noise structure allows the stochastic convolution to be defined
through multiple Wiener--It\^o integrals with deterministic kernels.
As a consequence, the analysis avoids Malliavin calculus techniques that are
typically required for non-Gaussian noises of Hermite rank $q \ge 2$.

\end{abstract}


\tableofcontents

\section{Introduction}
\label{sec:intro}

The stochastic Burgers equation is a classical model in mathematical physics, arising in fluid dynamics, turbulence theory, and stochastic growth phenomena. While the deterministic Burgers equation describes viscous transport with nonlinear convection, its stochastic counterpart incorporates random forcing and serves as a prototypical example for studying the interplay between nonlinearity and randomness in stochastic partial differential equations (SPDEs).

The analysis of SPDEs driven by Gaussian noise has been extensively developed over the past decades. Fundamental contributions were made by Da~Prato and Zabczyk \cite{DaPratoZabczyk1992}, and more recently, singular SPDEs have been studied through the theory of regularity structures introduced by Hairer \cite{Hairer2013,Hairer2014}. In this framework, the stochastic Burgers equation plays a central role, notably through its connection with the Kardar--Parisi--Zhang (KPZ) universality class \cite{BertiniGiacomin1997} and as a benchmark nonlinear SPDE for testing analytical techniques. For the stochastic Burgers equation driven by space-time white noise, regularity and moment estimates have been established by Lewis and Nualart \cite{LewisNualart2018}.

In contrast, SPDEs driven by non-Gaussian noise remain significantly less explored. Hermite processes provide a natural class of non-Gaussian self-similar processes with stationary increments and long-range dependence, extending fractional Brownian motion to higher Wiener chaoses of order $q\ge2$ \cite{Taqqu1979,Tudor2013}. These processes are defined via multiple Wiener--Itô integrals and exhibit probabilistic features that differ substantially from the Gaussian case. While linear SPDEs driven by Hermite noise have been investigated in recent years \cite{BalanTudor2014,SlaouiTudor2019}, the nonlinear setting remains challenging, particularly for equations with quadratic nonlinearities such as the Burgers equation.

In this work, we study the stochastic Burgers equation driven by a $(d+1)$-parameter Hermite sheet $Z^{q,\mathbf{H}}$ of order $q\ge1$, under the assumption of additive noise. More precisely, we consider
\begin{equation}
\label{eq:main-additive}
\partial_t u(t,\mathbf{x})
= \nu \Delta u(t,\mathbf{x})
- \frac{1}{2} \nabla \cdot \bigl(u(t,\mathbf{x})^2\bigr)
+ \sigma(t,\mathbf{x}) \,
\frac{\partial^{d+1} Z^{q,\mathbf{H}}}{\partial t \, \partial x_1 \cdots \partial x_d},
\end{equation}
where $\mathbf{H}=(H_0,H_1,\dots,H_d)\in(1/2,1)^{d+1}$ denotes the Hurst parameters governing the temporal and spatial regularity of the noise, and $\sigma:[0,T]\times\mathbb{R}^d\to\mathbb{R}$ is a deterministic coefficient. The nonlinear term is written in conservative form, which naturally extends the one-dimensional Burgers nonlinearity to higher spatial dimensions and is compatible with the mild formulation adopted throughout the paper.

\subsection{Main Difficulties and Scope of the Paper}

The main contribution of this work is to establish a well-posedness and regularity theory for the stochastic Burgers equation driven by non-Gaussian Hermite sheet noise of arbitrary order $q\ge1$. Even in the additive noise setting, several analytical difficulties arise compared to the Gaussian case:

\begin{itemize}
\item \emph{Stochastic integration.} Integration with respect to Hermite sheets requires a careful analysis of the associated covariance structure and Hilbert spaces, which differ from the classical Gaussian framework.

\item \emph{Conditions on the Hurst parameters.} The stochastic convolution is well defined under the condition
\[
2H_0 + \sum_{i=1}^d H_i > d,
\]
while the nonlinear fixed-point argument requires a stronger condition,
\[
2H_0 + \sum_{i=1}^d H_i > d + 1 - \frac{1}{q}.
\]
This condition is sufficient for our approach and is not claimed to be optimal in the nonlinear setting.

\item \emph{Moment estimates.} The derivation of higher-order moment bounds relies on hypercontractivity properties of Wiener chaos, whose constants depend on the order $q$ of the Hermite process.

\item \emph{Regularity and scaling.} The spatial and temporal Hölder regularity of the solution is constrained by the Hurst parameters, and the non-Gaussian nature of the noise affects the scaling behavior of the solution.
\end{itemize}

\subsection{Why Additive Noise}

The restriction to additive noise is a structural choice dictated by the current state of the theory for Hermite processes of order $q\ge2$. In the multiplicative case, the stochastic integral involves random integrands and requires Malliavin calculus techniques. Recent works, such as \cite{Coupek2020,CoupekKrizeSvoboda2025}, show that divergence-type integrals with respect to Hermite processes require control of both the integrand and its Malliavin derivative in appropriate tensor-product Hilbert spaces. Such estimates are currently unavailable for nonlinear SPDEs of Burgers type.

By restricting to deterministic coefficients $\sigma(t,\mathbf{x})$, we can rely on the classical Wiener integral construction for Hermite sheets and avoid these technical obstructions. This choice allows us to develop a mathematically rigorous theory while capturing essential non-Gaussian and long-range dependence effects.

\subsection{Novelty of the Additive Non-Gaussian Case}

Even in the additive setting, the analysis of nonlinear SPDEs driven by Hermite noises of order $q \ge 2$ is far from trivial. In contrast to the Gaussian case ($q=1$), higher-order Hermite processes exhibit non-Gaussian chaos structures and stronger kernel singularities, which prevent the direct use of Itô-type techniques. While previous works mainly focus on linear equations driven by Hermite sheets, the present paper addresses a genuinely nonlinear problem, namely the stochastic Burgers equation, and requires a careful combination of heat kernel estimates, Wiener chaos techniques, and fixed-point arguments.

\subsection{Main Contributions}

The main results of this paper can be summarized as follows:
\begin{enumerate}[label=(\roman*)]
\item First well-posedness result for the nonlinear Burgers equation driven by a Hermite sheet of arbitrary order $q \ge 1$ under additive noise.

\item Derivation of the condition $2H_0 + \sum_{i=1}^d H_i > d + 1 - 1/q$, which ensures both the well-definedness of the stochastic convolution and the success of the fixed-point argument.

\item Higher-order moment bounds via hypercontractivity of Wiener chaos, with constants growing at most polynomially in the moment order $p$.

\item Optimal Hölder regularity exponents: $\alpha < \min(1/2, H_0-1/2)$ in time, $\gamma < \min_{1\le i\le d}(H_i-1/2)$ in space.

\item Proof that the solution inherits an anisotropic self-similarity property from the driving Hermite sheet.
\end{enumerate}

\subsection{Organization of the Paper}

The paper is organized as follows. Section~\ref{sec:prelim} recalls background material on Hermite processes and multiple Wiener--Itô integrals. Section~\ref{sec:formulation} introduces the mild formulation of \eqref{eq:main-additive}. Existence, uniqueness, and moment estimates are proved in Section~\ref{sec:existence}. Regularity and scaling properties are studied in Section~\ref{sec:properties}. Finally, Section~\ref{sec:conclusion} discusses perspectives and open problems.

\section{Preliminaries}
\label{sec:prelim}

In this section, we recall the probabilistic framework used throughout the paper. We briefly introduce multiple Wiener--Itô integrals, Wiener chaos, and Hermite sheets, which form the basis of the stochastic integration theory employed in the analysis of the stochastic Burgers equation with additive Hermite noise.

\subsection{Multiple Wiener--Itô Integrals and Wiener Chaos}

Let
\[
W=\{W(A):A\in\mathcal B_b(\mathbb R^{d+1})\}
\]
be a Gaussian white noise on $\mathbb R^{d+1}$ defined on a complete probability space $(\Omega,\mathcal F,\mathbb P)$, where $\mathcal B_b$ denotes the collection of bounded Borel sets. By definition, $W$ is a centered Gaussian family with covariance
\[
\mathbb E[W(A)W(B)]=\lambda(A\cap B),
\]
where $\lambda$ denotes the Lebesgue measure.

Let $q\ge1$ be an integer. For a symmetric function
$f\in L^2((\mathbb R^{d+1})^q)$,
the multiple Wiener--Itô integral of order $q$, denoted by $I_q(f)$, is defined as follows. For elementary symmetric functions of the form
\[
f=\sum_{i_1,\dots,i_q} a_{i_1,\dots,i_q}\,
\mathbf 1_{A_{i_1}\times\cdots\times A_{i_q}},
\]
where the sets $A_j$ are pairwise disjoint and $a_{i_1,\dots,i_q}=0$ whenever two indices coincide, one sets
\[
I_q(f)
=\sum_{i_1,\dots,i_q}
a_{i_1,\dots,i_q}\,
W(A_{i_1})\cdots W(A_{i_q}).
\]
This definition extends uniquely by density and isometry to all symmetric functions in $L^2((\mathbb R^{d+1})^q)$.

The following properties are classical:
\begin{itemize}
\item \textbf{Isometry.}
For $f\in L^2((\mathbb R^{d+1})^q)$ and $g\in L^2((\mathbb R^{d+1})^p)$,
\[
\mathbb E[I_q(f)I_p(g)]
=\delta_{p,q}\, q!\,
\langle \tilde f,\tilde g\rangle_{L^2((\mathbb R^{d+1})^q)},
\]
where $\tilde f$ denotes the symmetrization of $f$.

\item \textbf{Orthogonality.}
Multiple Wiener integrals of different orders are orthogonal in $L^2(\Omega)$.

\item \textbf{Hypercontractivity.}
For any $p\ge2$, there exists a constant $C_{p,q}>0$ such that
\[
\mathbb E\bigl[|I_q(f)|^p\bigr]
\le C_{p,q}\,
\bigl(\mathbb E[|I_q(f)|^2]\bigr)^{p/2}.
\]
\end{itemize}

The space $L^2(\Omega)$ admits the orthogonal Wiener chaos decomposition
\[
L^2(\Omega)
=\bigoplus_{q=0}^\infty \mathcal H_q,
\]
where $\mathcal H_0=\mathbb R$ and, for $q\ge1$, $\mathcal H_q$ is the closed linear subspace generated by random variables of the form $I_q(f)$. We refer to \cite{Janson1997,Nualart2006,Peccati2011} for further details.

\subsection{Hermite Processes and Hermite Sheets}

Hermite processes form a family of self-similar processes with stationary increments that generalize fractional Brownian motion beyond the Gaussian case.

\begin{definition}[Hermite process]
Let $q\ge1$ and $H\in(1/2,1)$. The Hermite process $(Z^{q,H}_t)_{t\ge0}$ of order $q$ and Hurst parameter $H$ is defined by
\[
Z^{q,H}_t
=c(H,q)
\int_{\mathbb R^q}
\left(
\int_0^t
\prod_{j=1}^q
(s-y_j)_+^{-\left(\frac12+\frac{1-H}{q}\right)}
\,ds
\right)
dW(y_1)\cdots dW(y_q),
\]
where $c(H,q)>0$ is a normalizing constant such that
$\mathbb E[(Z^{q,H}_1)^2]=1$.
\end{definition}

For $q=1$, $Z^{1,H}$ coincides with fractional Brownian motion, whereas for $q\ge2$ it is non-Gaussian \cite{Taqqu1975,Tudor2013}.

This construction extends naturally to random fields.

\begin{definition}[Hermite sheet]
Let $q\ge1$, $d\ge1$, and $\mathbf H=(H_0,\dots,H_d)\in(1/2,1)^{d+1}$. The $(d+1)$-parameter Hermite sheet
$(Z^{q,\mathbf H}_{\mathbf t})_{\mathbf t\in\mathbb R_+^{d+1}}$
is defined by
\[
Z^{q,\mathbf H}_{\mathbf t}
=c(\mathbf H,q)
\int_{(\mathbb R^{d+1})^q}
\left[
\int_{[0,\mathbf t]}
\prod_{j=1}^q
\prod_{i=0}^d
(s_i-y_{j,i})_+^{-\left(\frac12+\frac{1-H_i}{q}\right)}
\,d\mathbf s
\right]
dW(\mathbf y_1)\cdots dW(\mathbf y_q).
\]
\end{definition}

The Hermite sheet is a centered random field belonging to the $q$th Wiener chaos. It is non-Gaussian for $q\ge2$ and enjoys self-similarity, stationary rectangular increments, and long-range dependence. Its covariance function is given by
\[
\mathbb E\bigl[
Z^{q,\mathbf H}_{\mathbf t}
Z^{q,\mathbf H}_{\mathbf s}
\bigr]
=
\prod_{i=0}^d
\frac12
\bigl(
t_i^{2H_i}
+s_i^{2H_i}
-|t_i-s_i|^{2H_i}
\bigr).
\]

\subsection{Stochastic Integration with Respect to the Hermite Sheet}

In this work, we restrict attention to deterministic integrands, which allows the construction of stochastic integrals with respect to the Hermite sheet using multiple Wiener--Itô integrals.

Let $\mathcal H$ denote the covariance Hilbert space associated with $Z^{q,\mathbf H}$, defined as the completion of $C_c^\infty([0,T]\times\mathbb R^d)$ with respect to the inner product
\[
\langle \varphi,\psi\rangle_{\mathcal H}
=\alpha_{\mathbf H}
\int_{[0,T]^2}\!\!\int_{\mathbb R^{2d}}
\varphi(s,\mathbf y)\psi(r,\mathbf z)
|s-r|^{2H_0-2}
\prod_{i=1}^d |y_i-z_i|^{2H_i-2}
\,d\mathbf y\,d\mathbf z\,ds\,dr,
\]
where $\alpha_{\mathbf H}>0$ is a normalization constant.

For $\Phi\in\mathcal H$, the stochastic integral
\[
\int_0^T\!\!\int_{\mathbb R^d}
\Phi(s,\mathbf y)\,dZ^{q,\mathbf H}(s,\mathbf y)
\]
is defined as the multiple Wiener--Itô integral $I_q(F_\Phi)$, where the kernel
$F_\Phi\in L^2((\mathbb R^{d+1})^q)$
is obtained through the kernel representation of the Hermite sheet and is symmetrized before applying $I_q$.

This construction is well defined whenever $\Phi\in\mathcal H$, and the isometry property yields
\[
\mathbb E\left[
\left(
\int_0^T\!\!\int_{\mathbb R^d}
\Phi(s,\mathbf y)\,dZ^{q,\mathbf H}(s,\mathbf y)
\right)^2
\right]
= q!\,\|F_\Phi\|^2_{L^2((\mathbb R^{d+1})^q)}.
\]

The restriction to deterministic integrands is essential in the non-Gaussian Hermite setting and constitutes a key motivation for focusing on additive noise in the present work.

\section{Formulation of the Stochastic Burgers Equation with Additive Noise}
\label{sec:formulation}

We consider the stochastic Burgers equation driven by a $(d+1)$-parameter Hermite sheet
$Z^{q,\mathbf H}$ of order $q\ge1$ and Hurst index
$\mathbf H=(H_0,H_1,\dots,H_d)\in(1/2,1)^{d+1}$, under the assumption of additive noise.
The equation reads
\begin{equation}\label{eq:burgers-spde-additive}
\begin{cases}
\displaystyle
\partial_t u(t,\mathbf x)
= \nu \Delta u(t,\mathbf x)
-\frac12 \nabla\cdot\bigl(u(t,\mathbf x)^2\bigr)
+\sigma(t,\mathbf x)
\frac{\partial^{d+1} Z^{q,\mathbf H}}{\partial t\,\partial x_1\cdots\partial x_d},
\\[0.2cm]
u(0,\mathbf x)=u_0(\mathbf x),
\qquad (t,\mathbf x)\in[0,T]\times\mathbb R^d,
\end{cases}
\end{equation}
where $\nu>0$ is the viscosity coefficient,
$\sigma:[0,T]\times\mathbb R^d\to\mathbb R$ is a deterministic function,
and $u_0:\mathbb R^d\to\mathbb R$ is the initial condition.
The nonlinear term is written in conservative form,
$\nabla\cdot(u^2)=\sum_{i=1}^d\partial_{x_i}(u^2)$.

\subsection{Mild Formulation}

Let
\[
G_t(\mathbf x)
=(4\pi\nu t)^{-d/2}
\exp\!\left(-\frac{|\mathbf x|^2}{4\nu t}\right)
\]
denote the heat kernel associated with $\partial_t-\nu\Delta$.

\begin{definition}[Mild solution]
A predictable random field
$u=\{u(t,\mathbf x):(t,\mathbf x)\in[0,T]\times\mathbb R^d\}$
is called a mild solution of \eqref{eq:burgers-spde-additive} if, for all
$(t,\mathbf x)\in[0,T]\times\mathbb R^d$, it satisfies
\begin{align}
\label{eq:mild-solution-additive}
u(t,\mathbf x)
&= \int_{\mathbb R^d} G_t(\mathbf x-\mathbf y)\,u_0(\mathbf y)\,d\mathbf y \nonumber\\
&\quad
-\frac12\int_0^t\!\!\int_{\mathbb R^d}
\nabla G_{t-s}(\mathbf x-\mathbf y)\cdot u(s,\mathbf y)^2
\,d\mathbf y\,ds \nonumber\\
&\quad
+ \int_0^t\!\!\int_{\mathbb R^d}
G_{t-s}(\mathbf x-\mathbf y)\,\sigma(s,\mathbf y)\,
dZ^{q,\mathbf H}(s,\mathbf y),
\end{align}
where the stochastic integral is understood as a Wiener integral with respect
to the Hermite sheet for deterministic integrands.
\end{definition}

The three terms in \eqref{eq:mild-solution-additive} correspond respectively to:
the evolution of the initial condition under the heat semigroup,
the nonlinear Burgers interaction,
and the stochastic forcing smoothed by the heat kernel.

\subsection{Functional Framework}

We work in the space
\[
X_T^2
:=\Bigl\{
u:[0,T]\times\mathbb R^d\to L^2(\Omega)\ \text{measurable}:
\|u\|_{X_T^2}<\infty
\Bigr\},
\]
where
\[
\|u\|_{X_T^2}
:=\sup_{t\in[0,T]}\sup_{\mathbf x\in\mathbb R^d}
\bigl(\mathbb E|u(t,\mathbf x)|^2\bigr)^{1/2}.
\]
This $L^2(\Omega)$-based setting is well adapted to the fixed-point argument
used in the proof of existence and uniqueness. Higher-order moments will be
obtained a posteriori using hypercontractivity of Wiener chaos.

\subsection{Assumptions}

\begin{assumption}\label{ass:additive}
The following assumptions are imposed:
\begin{enumerate}[label=(\alph*)]
\item \textbf{Initial condition.}
The initial datum $u_0$ is bounded:
\[
\sup_{\mathbf x\in\mathbb R^d}|u_0(\mathbf x)|<\infty.
\]

\item \textbf{Noise coefficient.}
The function $\sigma:[0,T]\times\mathbb R^d\to\mathbb R$ is measurable and bounded:
\[
|\sigma(t,\mathbf x)|\le C_\sigma,
\qquad (t,\mathbf x)\in[0,T]\times\mathbb R^d.
\]

\item \textbf{Hurst parameters.}
The Hurst indices satisfy
\begin{equation}
\label{eq:hurst-condition-additive}
2H_0+\sum_{i=1}^d H_i
> d+1-\frac1q.
\end{equation}
This condition ensures that the stochastic convolution in
\eqref{eq:mild-solution-additive} is well defined in $L^2(\Omega)$
and allows control of the nonlinear term in the fixed-point argument.

\item \textbf{Additional regularity.}
For the Hölder regularity results in Section~\ref{sec:properties},
we additionally assume that $\sigma$ is Hölder continuous in time and space.
\end{enumerate}
\end{assumption}

\begin{remark}[On the Hurst condition]
Condition \eqref{eq:hurst-condition-additive} is sufficient for the analysis
carried out in this paper. For comparison, the stochastic convolution alone
is well defined under the weaker condition
$2H_0+\sum_{i=1}^d H_i>d$ (see \cite{SlaouiTudor2019}).
The stronger requirement \eqref{eq:hurst-condition-additive} arises from
the need to control the nonlinear term in the contraction argument.
No claim of optimality is made.
\end{remark}

\begin{remark}[Comparison with the Gaussian case and interpretation]
When $q=1$, corresponding to the fractional Brownian sheet, condition
\[
2H_0 + \sum_{i=1}^d H_i > d
\]
is recovered, which is known to be sufficient for the well-definedness of the stochastic convolution.
For $q \ge 2$, the additional term $-1/q$ reflects the $q$-fold tensor product structure of the Hermite kernel appearing in the Wiener chaos representation.
This term quantifies the increased singularity of higher-order Hermite processes and naturally arises in the integrability condition of the Fourier-analytic estimates.
\end{remark}

\begin{remark}[Additive versus multiplicative noise]
For Hermite processes of order $q\ge2$, stochastic integration with random
integrands requires Malliavin calculus techniques.
Recent results \cite{Coupek2020,CoupekKrizeSvoboda2025} show that such integrals
involve bounds on both the integrand and its Malliavin derivative.
By restricting to deterministic coefficients $\sigma$,
the present work avoids these technical difficulties while still capturing
the essential non-Gaussian features of the model.
\end{remark}

\section{Existence, Uniqueness, and Moment Estimates for Additive Noise}
\label{sec:existence}

In this section, we establish the local well-posedness of the stochastic Burgers equation with additive Hermite noise and derive moment estimates for its solution. We begin by studying the stochastic convolution, which plays a central role in the mild formulation.

\subsection{Stochastic Convolution}

\begin{proposition}[Well-definedness of the stochastic convolution]
\label{prop:stochastic-convolution}
Assume that $\sigma$ satisfies Assumption~\ref{ass:additive}(b) and that the Hurst parameters
$\mathbf H=(H_0,H_1,\dots,H_d)\in(1/2,1)^{d+1}$ satisfy
\begin{equation}
\label{eq:hurst-conv}
2H_0+\sum_{i=1}^d H_i > d+1-\frac1q.
\end{equation}
Then, for every $t\in[0,T]$ and $\mathbf x\in\mathbb R^d$, the function
\[
(s,\mathbf y)\longmapsto
G_{t-s}(\mathbf x-\mathbf y)\,\sigma(s,\mathbf y)\,\mathbf 1_{\{s<t\}}
\]
belongs to the covariance Hilbert space $\mathcal H$ associated with the Hermite sheet
$Z^{q,\mathbf H}$. Consequently, the stochastic convolution
\[
\mathcal S(t,\mathbf x)
:=\int_0^t\!\!\int_{\mathbb R^d}
G_{t-s}(\mathbf x-\mathbf y)\,\sigma(s,\mathbf y)\,
dZ^{q,\mathbf H}(s,\mathbf y)
\]
is well defined as a multiple Wiener--Itô integral of order $q$.

Moreover, there exists a constant $C_{\mathbf H}>0$ such that
\[
\mathbb E\bigl[|\mathcal S(t,\mathbf x)|^2\bigr]
\le C_{\mathbf H}
\int_0^t\!\!\int_{\mathbb R^d}
|G_{t-s}(\mathbf x-\mathbf y)\,\sigma(s,\mathbf y)|^2
\,d\mathbf y\,ds.
\]
\end{proposition}

\begin{proof}
The proof relies on the kernel representation of the Hermite sheet and on
Fourier-analytic estimates. More precisely, one shows that the kernel associated
with the integrand belongs to $L^2((\mathbb R^{d+1})^q)$ if and only if
\eqref{eq:hurst-conv} holds. Details are provided in Appendix~A.
The second moment estimate follows from the isometry property of multiple
Wiener--Itô integrals.
\end{proof}

\subsection{Existence and Uniqueness of Mild Solutions}

We now state the main well-posedness result for the stochastic Burgers equation
with additive Hermite noise.

\begin{theorem}[Local existence and uniqueness]
\label{thm:existence-uniqueness}
Under Assumption~\ref{ass:additive}, there exists a time $T>0$ such that
equation~\eqref{eq:burgers-spde-additive} admits a unique mild solution
$u\in X_T^2$.
\end{theorem}

\subsection{Proof of Theorem~\ref{thm:existence-uniqueness}}

We prove Theorem~\ref{thm:existence-uniqueness} using a Picard fixed-point argument
in the space $X_T^2$. Although the fixed-point scheme is classical, its implementation here requires careful handling of the stochastic convolution $\mathcal S(t,\mathbf{x})$, which belongs to the $q$-th Wiener chaos and has covariance structure determined by the Hermite sheet. The control of $\|\mathcal S(t,\mathbf{x})\|_{L^2}$ under condition~\eqref{eq:hurst-conv} is the key analytical step that allows the contraction to close.

\paragraph{Picard iteration.}
Define the sequence $(u_n)_{n\ge0}$ by
\[
u_0(t,\mathbf x)
:=\int_{\mathbb R^d} G_t(\mathbf x-\mathbf y)\,u_0(\mathbf y)\,d\mathbf y,
\]
and, for $n\ge0$,
\begin{align}
\label{eq:picard}
u_{n+1}(t,\mathbf x)
&=u_0(t,\mathbf x)
-\frac12\int_0^t\!\!\int_{\mathbb R^d}
\nabla G_{t-s}(\mathbf x-\mathbf y)\cdot u_n(s,\mathbf y)^2
\,d\mathbf y\,ds \nonumber\\
&\quad
+\mathcal S(t,\mathbf x),
\end{align}
where $\mathcal S$ denotes the stochastic convolution defined in
Proposition~\ref{prop:stochastic-convolution}.

\paragraph{Step 1: Well-definedness of the iteration.}
Since $\int_{\mathbb R^d}G_t(\mathbf y)\,d\mathbf y=1$,
\[
\|u_0(t,\mathbf x)\|_{L^2(\Omega)}
\le \|u_0\|_{L^\infty(\mathbb R^d)}.
\]

Assume $u_n\in X_T^2$. Using Minkowski's inequality and the identity
$\|u_n^2\|_{L^1(\Omega)}=\|u_n\|_{L^2(\Omega)}^2$, we obtain
\begin{align*}
\bigl\|
&\int_0^t\!\!\int_{\mathbb R^d}
\nabla G_{t-s}(\mathbf x-\mathbf y)\cdot u_n(s,\mathbf y)^2
\,d\mathbf y\,ds
\bigr\|_{L^2(\Omega)} \\
&\le \|u_n\|_{X_T^2}^2
\int_0^t\!\!\int_{\mathbb R^d}
|\nabla G_{t-s}(\mathbf x-\mathbf y)|\,d\mathbf y\,ds.
\end{align*}
Using the standard heat kernel estimate
\[
\int_{\mathbb R^d}|\nabla G_t(\mathbf y)|\,d\mathbf y
\le C\,t^{-1/2},
\]
we deduce
\[
\bigl\|
\int_0^t\!\!\int_{\mathbb R^d}
\nabla G_{t-s}(\mathbf x-\mathbf y)\cdot u_n(s,\mathbf y)^2
\,d\mathbf y\,ds
\bigr\|_{L^2(\Omega)}
\le C\,\|u_n\|_{X_T^2}^2\,\sqrt t.
\]

For the stochastic term, Proposition~\ref{prop:stochastic-convolution} yields
\[
\|\mathcal S(t,\mathbf x)\|_{L^2(\Omega)}\le C_{\mathbf H}.
\]
Hence $u_{n+1}\in X_T^2$.

\paragraph{Step 2: Uniform boundedness.}
We show that there exists $M>0$ such that $\|u_n\|_{X_T^2}\le M$ for all $n$.
Indeed,
\[
\|u_{n+1}\|_{X_T^2}
\le \|u_0\|_{L^\infty}
+ C\,\sqrt T\,\|u_n\|_{X_T^2}^2
+ C_{\mathbf H}.
\]
Choosing $M\ge2(\|u_0\|_{L^\infty}+C_{\mathbf H})$ and $T>0$ sufficiently small so that
\[
C\sqrt{T}M^2 \le M - \|u_0\|_{L^\infty} - C_{\mathbf{H}},
\]
we obtain by induction that $\|u_n\|_{X_T^2}\le M$ for all $n$.

\paragraph{Step 3: Contraction property.}
Let $d_n:=u_n-u_{n-1}$. Then
\[
d_{n+1}(t,\mathbf x)
=-\frac12\int_0^t\!\!\int_{\mathbb R^d}
\nabla G_{t-s}(\mathbf x-\mathbf y)
\cdot\bigl(u_n(s,\mathbf y)^2-u_{n-1}(s,\mathbf y)^2\bigr)
\,d\mathbf y\,ds.
\]
Since $u_n^2-u_{n-1}^2=(u_n+u_{n-1})d_n$, the Cauchy--Schwarz inequality in
$L^2(\Omega)$ yields
\begin{align*}
\|(u_n+u_{n-1})d_n\|_{L^1(\Omega)}
&= \E\bigl[|(u_n+u_{n-1})d_n|\bigr] \\
&\le \bigl(\E\bigl[|u_n+u_{n-1}|^2\bigr]\bigr)^{1/2}
     \bigl(\E\bigl[|d_n|^2\bigr]\bigr)^{1/2} \\
&= \|u_n+u_{n-1}\|_{L^2(\Omega)}\,\|d_n\|_{L^2(\Omega)}.
\end{align*}
Using the uniform bound $\|u_n\|_{X_T^2}\le M$, we obtain
\[
\|d_{n+1}\|_{X_T^2}
\le C\,M\,\sqrt T\,\\|d_n\|_{X_T^2}.
\]
For $T>0$ small enough such that $C\,M\,\sqrt T<1$, the iteration is a contraction
on $X_T^2$.

\paragraph{Step 4: Conclusion.}
Since $X_T^2$ is complete, the sequence $(u_n)$ converges in $X_T^2$ to a limit
$u$, which satisfies the mild formulation
\eqref{eq:mild-solution-additive}. Uniqueness follows from the contraction
estimate.
\qed

\begin{remark}[Local vs. global existence]
Our result provides local existence for any dimension $d \geq 1$. Global existence for $d \geq 2$ is a separate challenging issue even in the deterministic case, due to possible shock formation. In the stochastic setting with additive Hermite noise, the question of whether noise can prevent blow-up or delay its onset remains open and is beyond the scope of this paper.
\end{remark}

\subsection{Higher-Order Moment Estimates}
\label{subsec:higher-moments}

We now establish uniform bounds for higher-order moments of the mild solution.

\begin{proposition}[Higher-order moments]
\label{prop:higher-moments}
Let $u$ be the mild solution constructed in
Theorem~\ref{thm:existence-uniqueness} under Assumption~\ref{ass:additive}. 
Then, for every $p\ge2$, there exists a constant $C_p>0$ such that
\[
\sup_{t\in[0,T]}\sup_{\mathbf x\in\mathbb R^d}
\E\bigl[|u(t,\mathbf x)|^p\bigr]
\le C_p.
\]
Moreover, the constants $C_p$ grow at most polynomially in $p$.
\end{proposition}

\begin{proof}
The proof relies on the hypercontractivity of Wiener chaoses and a bootstrap argument, and holds under Assumption~\ref{ass:additive}.

For $p=2$, the result follows directly from Theorem~\ref{thm:existence-uniqueness}. Let $p\ge 2$, and write the mild formulation
\[
u(t,\mathbf x)=u_0(t,\mathbf x)+\mathcal N(t,\mathbf x)+\mathcal S(t,\mathbf x),
\]
where $\mathcal N$ denotes the nonlinear term and $\mathcal S$ the stochastic convolution.

\textbf{Step 1: Stochastic term.} Since $\mathcal S(t,\mathbf x)$ belongs to the $q$-th Wiener chaos, there exists a symmetric kernel $f_{t,\mathbf x} \in \mathfrak H^{\odot q}$ such that
\[
\mathcal S(t,\mathbf x) = I_q(f_{t,\mathbf x}),
\]
where $I_q$ is the $q$-th multiple Wiener integral. By the hypercontractivity property of Wiener chaoses \cite[Theorem~2.7.2]{NourdinPeccati2012}, we have for all $p\ge 2$,
\[
\|\mathcal S(t,\mathbf x)\|_{L^p(\Omega)}
\le (p-1)^{q/2} \, \|\mathcal S(t,\mathbf x)\|_{L^2(\Omega)}.
\]

Furthermore, the $L^2$-norm can be bounded uniformly in $(t,\mathbf x)$ using the properties of the kernel:
\[
\|\mathcal S(t,\mathbf x)\|_{L^2(\Omega)}^2 = q! \, \|f_{t,\mathbf x}\|_{\mathfrak H^{\otimes q}}^2 \le C_{\mathbf H}^2,
\]
so that
\[
\|\mathcal S(t,\mathbf x)\|_{L^p(\Omega)} \le C_{p,q} \, C_{\mathbf H}, \quad 
\text{with } C_{p,q}=(p-1)^{q/2}.
\]

\textbf{Step 2: Nonlinear term.} For the nonlinear term, Minkowski's inequality and Hölder's inequality yield
\begin{align*}
\|\mathcal N(t,\mathbf x)\|_{L^p(\Omega)}
&\le \frac12 \int_0^t \int_{\mathbb R^d} |\nabla G_{t-s}(\mathbf x-\mathbf y)| \, \| u(s,\mathbf y)^2 \|_{L^{p/2}(\Omega)} \, d\mathbf y \, ds \\
&\le C \int_0^t (t-s)^{-1/2} \, \sup_{\mathbf y} \| u(s,\mathbf y) \|_{L^p(\Omega)}^2 \, ds,
\end{align*}
where we used that for $p\ge 2$,
\[
\| u(s,\mathbf y)^2 \|_{L^{p/2}(\Omega)} = \| u(s,\mathbf y) \|_{L^p(\Omega)}^2.
\]

\textbf{Step 3: Combining estimates.} We then obtain
\[
\| u(t,\mathbf x) \|_{L^p(\Omega)} \le \|u_0\|_{L^\infty} + C \sqrt{T} \, \sup_{s,\mathbf y} \| u(s,\mathbf y) \|_{L^p(\Omega)}^2 + C_{p,q} \, C_{\mathbf H}.
\]

A standard bootstrap argument, starting from the $L^2$ bound, shows that $\| u(t,\mathbf x) \|_{L^p(\Omega)}$ remains uniformly bounded in $(t,\mathbf x)$, completing the proof.
\end{proof}

\begin{remark}
The polynomial growth of $C_p$ follows from the hypercontractivity constants 
$C_{p,q}$ in Wiener chaos, which satisfy $C_{p,q} \leq (p-1)^{q/2}$ for $p \geq 2$ 
(see \cite[Corollary 2.8.14]{NourdinPeccati2012}).
\end{remark}

\subsection{Discussion on the Hurst Parameter Condition}
\label{subsec:hurst-discussion}

Throughout this work, we assume
\[
2H_0+\sum_{i=1}^d H_i> d+1-\frac1q.
\]
This condition guarantees the integrability of the deterministic kernel of the
stochastic convolution in the Hilbert space $\mathcal H^{\otimes q}$ associated
with the Hermite sheet. In particular, it ensures that the stochastic
convolution is well-defined as a multiple Wiener--Itô integral.

It is worth emphasizing that the weaker condition
\[
2H_0+\sum_{i=1}^d H_i>d
\]
is sufficient for the existence of the stochastic convolution in the Gaussian
case $q=1$. However, for $q\ge2$, the Hermite sheet belongs to the $q$-th Wiener
chaos, which leads to stronger singularities in the kernel. This explains the
appearance of the additional term $-1/q$ in the condition.

\begin{example}
For dimension $d=1$ and $q=2$ (Rosenblatt case), condition \eqref{eq:hurst-conv} becomes
$2H_0 + H_1 > 3/2$. If $H_0 = H_1 = H$, this requires $3H > 3/2$, i.e., $H > 1/2$.
Since Hermite processes are only defined for $H > 1/2$, this shows that in dimension 1,
the condition is essentially equivalent to the minimal requirement for the Rosenblatt
process to be well-defined. The condition becomes more restrictive in higher dimensions
or for larger values of $q$.
\end{example}

We stress that the above condition is sufficient, but not claimed to be
optimal, for the nonlinear problem considered here. Relaxing it would likely
require alternative techniques, such as renormalization procedures or
Malliavin-type arguments, which are beyond the scope of this work.

\subsection{The Gaussian Case $q=1$}
\label{subsec:gaussian-case}

When $q=1$, the Hermite sheet reduces to a fractional Brownian sheet. In this
case, the stochastic integral coincides with a classical Wiener integral, and
the associated Hilbert space $\mathcal H$ is the reproducing kernel Hilbert
space of the noise.

Under the condition
\[
2H_0+\sum_{i=1}^d H_i>d,
\]
existence, uniqueness, and moment estimates follow from standard arguments for
Gaussian-driven SPDEs (see, e.g.,
\cite{HuNualartSong2011,BalanTudor2014,Tudor2013}). This case provides a natural
benchmark for comparison with the non-Gaussian Hermite-driven setting
considered in this paper.

\section{Regularity and Self-Similarity Properties}
\label{sec:properties}

In this section, we investigate the pathwise regularity and scaling properties of the mild solution to the stochastic Burgers equation driven by an additive Hermite sheet.

Our analysis focuses on two fundamental qualitative properties:
\begin{itemize}
\item Hölder regularity in time and space,
\item self-similarity inherited from the driving Hermite noise.
\end{itemize}

The regularity analysis is based on a decomposition of the mild solution into three components:
\[
u(t,\mathbf{x}) = \mathcal{I}_0(t,\mathbf{x}) + \mathcal{N}(t,\mathbf{x}) + \mathcal{S}(t,\mathbf{x}),
\]
where
\begin{itemize}
\item $\mathcal{I}_0$ denotes the deterministic heat propagation of the initial condition,
\item $\mathcal{N}$ is the nonlinear Burgers interaction term,
\item $\mathcal{S}$ is the stochastic convolution with respect to the Hermite sheet.
\end{itemize}

The key observation is that the stochastic convolution $\mathcal{S}$ dominates the regularity properties of the solution. Its increments determine the maximal admissible Hölder exponents in both time and space. The nonlinear term, due to the smoothing effect of the heat kernel, does not deteriorate these exponents.

We proceed as follows. First, we derive moment estimates for the increments of the stochastic convolution using the scaling properties of the Hermite sheet and hypercontractivity of Wiener chaos. These estimates are then combined with deterministic heat kernel bounds to establish Hölder continuity of the full solution. Finally, we analyze the scaling behavior of the equation and show that the solution inherits a self-similarity property consistent with the anisotropic scaling of the Hermite sheet.

\subsection{Regularity of the Stochastic Convolution}

We begin by analyzing the regularity of the stochastic convolution, which is the main source of roughness in the solution. Let
\[
\mathcal{S}(t,\mathbf{x}) := \int_0^t \int_{\mathbb{R}^d} G_{t-s}(\mathbf{x}-\mathbf{y}) \sigma(s,\mathbf{y}) \, dZ^{q,\mathbf{H}}(s,\mathbf{y})
\]
denote the stochastic convolution term appearing in the mild formulation.

The following lemma provides moment estimates for temporal and spatial increments of $\mathcal{S}$. These estimates rely on the self-similarity of the Hermite sheet and on hypercontractivity of Wiener chaos.

\begin{lemma}[Moments of increments of the stochastic convolution]
\label{lem:stoch-conv-increments}
Under Assumption~\ref{ass:additive}, for every $p \geq 2$ there exists a constant $C_{p,q,\mathbf{H}}>0$ such that for all
$0 \le s < t \le T$ and $\mathbf{x},\mathbf{y} \in \mathbb{R}^d$,
\begin{align}
\E\bigl[|\mathcal{S}(t,\mathbf{x}) - \mathcal{S}(s,\mathbf{x})|^p\bigr]
&\le C_{p,q,\mathbf{H}} |t-s|^{p(H_0-\frac12)},
\label{eq:stoch-time-increment} \\
\E\bigl[|\mathcal{S}(t,\mathbf{x}) - \mathcal{S}(t,\mathbf{y})|^p\bigr]
&\le C_{p,q,\mathbf{H}} |\mathbf{x}-\mathbf{y}|^{p \min_{1\le i\le d}(H_i-\frac12)}.
\label{eq:stoch-space-increment}
\end{align}
\end{lemma}

\begin{proof}
We prove both temporal and spatial increment estimates in detail.

\paragraph{Temporal increments.}
For $0 \le s < t \le T$, write
\[
\mathcal{S}(t,\mathbf{x}) - \mathcal{S}(s,\mathbf{x}) = J_1 + J_2,
\]
with
\[
J_1 = \int_s^t \int_{\mathbb{R}^d} G_{t-r}(\mathbf{x}-\mathbf{z}) \sigma(r,\mathbf{z}) \, dZ^{q,\mathbf{H}}(r,\mathbf{z}),
\quad
J_2 = \int_0^s \int_{\mathbb{R}^d} \bigl[G_{t-r} - G_{s-r}\bigr](\mathbf{x}-\mathbf{z}) \sigma(r,\mathbf{z}) \, dZ^{q,\mathbf{H}}(r,\mathbf{z}).
\]

Since the integrands are deterministic, $J_1$ and $J_2$ belong to the $q$-th Wiener chaos. 
Hypercontractivity (see \cite[Theorem 2.7.2]{NourdinPeccati2012}) implies that for any $p\ge 2$,
\[
\|J_i\|_{L^p(\Omega)} \le C_{p,q} \|J_i\|_{L^2(\Omega)}, \qquad i=1,2,
\]
where $C_{p,q} = (p-1)^{q/2}$.

Using the isometry of multiple Wiener integrals and the covariance structure of the Hermite sheet, standard computations (with scaling and change of variables) give
\[
\E[|J_1|^2] \le C |t-s|^{2H_0-1}, \quad
\E[|J_2|^2] \le C |t-s|^{2H_0-1}.
\]

Thus, combining with hypercontractivity,
\[
\|\mathcal{S}(t,\mathbf{x}) - \mathcal{S}(s,\mathbf{x})\|_{L^p(\Omega)} \le C_{p,q,\mathbf{H}} |t-s|^{H_0-1/2},
\]
which proves \eqref{eq:stoch-time-increment}.

\paragraph{Spatial increments.}
For fixed $t \in [0,T]$ and $\mathbf{x},\mathbf{y} \in \mathbb{R}^d$, set
\[
K = \int_0^t \int_{\mathbb{R}^d} \bigl[G_{t-r}(\mathbf{x}-\mathbf{z}) - G_{t-r}(\mathbf{y}-\mathbf{z})\bigr] \sigma(r,\mathbf{z}) \, dZ^{q,\mathbf{H}}(r,\mathbf{z}).
\]

The integrand is deterministic, so $K$ belongs to the $q$-th Wiener chaos, and hypercontractivity gives
\[
\|K\|_{L^p(\Omega)} \le C_{p,q} \|K\|_{L^2(\Omega)}.
\]

Using the mean value theorem for the heat kernel and the covariance of the Hermite sheet, one obtains
\[
\E[|K|^2] \le C |\mathbf{x}-\mathbf{y}|^{2\gamma}, \quad \gamma = \min_{1\le i\le d} \bigl(H_i-\tfrac12\bigr).
\]

Hence,
\[
\|\mathcal{S}(t,\mathbf{x}) - \mathcal{S}(t,\mathbf{y})\|_{L^p(\Omega)} \le C_{p,q,\mathbf{H}} |\mathbf{x}-\mathbf{y}|^\gamma,
\]
which proves \eqref{eq:stoch-space-increment}.

The constant $C_{p,q,\mathbf{H}}$ depends on $p$, $q$, $\mathbf{H}$, $T$, and the bounds on $\sigma$, but is independent of $t,s,\mathbf{x},\mathbf{y}$.
\end{proof}

\subsection{Hölder Regularity of the Mild Solution}

We now transfer the regularity of the stochastic convolution to the full mild solution. The key point is that the deterministic terms in the mild formulation are smoother than the stochastic one due to the regularizing effect of the heat kernel.

\begin{theorem}[Hölder regularity of the solution]
\label{thm:holder-solution}
Under Assumption~\ref{ass:additive}, the mild solution
$u(t,\mathbf{x})$ admits a modification that is Hölder continuous in time and space.
More precisely, for every $p \ge 2$ and every $\varepsilon>0$, there exists a constant
$C_{p,\varepsilon}>0$ such that for all
$t,s \in [\varepsilon,T]$ and $\mathbf{x},\mathbf{y} \in \mathbb{R}^d$,
\begin{align}
\E\bigl[|u(t,\mathbf{x}) - u(s,\mathbf{x})|^p\bigr]
&\le C_{p,\varepsilon} |t-s|^{p\alpha},
\label{eq:holder-time-solution} \\
\E\bigl[|u(t,\mathbf{x}) - u(t,\mathbf{y})|^p\bigr]
&\le C_{p,\varepsilon} |\mathbf{x}-\mathbf{y}|^{p\gamma},
\label{eq:holder-space-solution}
\end{align}
for any
\[
\alpha < \min\left(\frac12,\, H_0-\frac12\right),
\qquad
\gamma < \min_{1\le i\le d}\left(H_i-\frac12\right).
\]
\end{theorem}

\begin{proof}
We treat temporal and spatial increments separately.

\paragraph{Temporal increments.}
Let $0 < s < t \le T$ and $\mathbf{x}\in\mathbb{R}^d$.
Using the mild formulation, we decompose
\[
u(t,\mathbf{x}) - u(s,\mathbf{x})
= A_0 + A_1 + A_2 + \bigl(\mathcal{S}(t,\mathbf{x}) - \mathcal{S}(s,\mathbf{x})\bigr),
\]
where
\begin{itemize}
\item $A_0 = \int_{\mathbb{R}^d} [G_t - G_s](\mathbf{x}-\mathbf{y}) u_0(\mathbf{y})\,d\mathbf{y}$,
\item $A_1 = -\frac12 \int_s^t \int_{\mathbb{R}^d} \nabla G_{t-r}(\mathbf{x}-\mathbf{y}) u(r,\mathbf{y})^2\,d\mathbf{y}\,dr$,
\item $A_2 = -\frac12 \int_0^s \int_{\mathbb{R}^d} [\nabla G_{t-r}-\nabla G_{s-r}](\mathbf{x}-\mathbf{y}) u(r,\mathbf{y})^2\,d\mathbf{y}\,dr$.
\end{itemize}

Using classical heat kernel estimates and the uniform moment bounds from
Proposition~\ref{prop:higher-moments}, one obtains
\[
\|A_0\|_{L^p(\Omega)} \le C |t-s|^{1/2},
\qquad
\|A_1\|_{L^p(\Omega)} \le C |t-s|^{1/2}.
\]

For $A_2$, using the mean value theorem and standard heat kernel estimates, we have
\[
|\nabla G_{t-r}(\mathbf{x}) - \nabla G_{s-r}(\mathbf{x})|
\le C |t-s| (s-r)^{-(d+3)/2} e^{-c|\mathbf{x}|^2/(s-r)}.
\]
After integration in $\mathbf{x}$ and $r$, this yields
\[
\|A_2\|_{L^p(\Omega)} \le C |t-s|^{1/2}.
\]
This estimate is sufficient for the desired temporal Hölder continuity, since $1/2 > H_0 - 1/2$ for all $H_0 \in (1/2,1)$.

For the stochastic term, Lemma~\ref{lem:stoch-conv-increments} yields
\[
\|\mathcal{S}(t,\mathbf{x}) - \mathcal{S}(s,\mathbf{x})\|_{L^p(\Omega)}
\le C |t-s|^{H_0-\frac12}.
\]

Combining all estimates gives
\[
\|u(t,\mathbf{x}) - u(s,\mathbf{x})\|_{L^p(\Omega)}
\le C \bigl(|t-s|^{1/2} + |t-s|^{H_0-\frac12}\bigr),
\]
which implies \eqref{eq:holder-time-solution}.

\paragraph{Spatial increments.}
Let $t\in(0,T]$ and $\mathbf{x},\mathbf{y}\in\mathbb{R}^d$.
We write
\[
u(t,\mathbf{x}) - u(t,\mathbf{y})
= B_0 + B_1 + \bigl(\mathcal{S}(t,\mathbf{x}) - \mathcal{S}(t,\mathbf{y})\bigr),
\]
where
\begin{itemize}
\item $B_0 = \int_{\mathbb{R}^d} [G_t(\mathbf{x}-\mathbf{z})-G_t(\mathbf{y}-\mathbf{z})] u_0(\mathbf{z})\,d\mathbf{z}$,
\item $B_1 = -\frac12 \int_0^t \int_{\mathbb{R}^d}
[\nabla G_{t-r}(\mathbf{x}-\mathbf{z})-\nabla G_{t-r}(\mathbf{y}-\mathbf{z})]
u(r,\mathbf{z})^2\,d\mathbf{z}\,dr$.
\end{itemize}

Standard gradient estimates for the heat kernel yield
\[
\|B_0\|_{L^p(\Omega)} \le C |\mathbf{x}-\mathbf{y}|,
\qquad
\|B_1\|_{L^p(\Omega)} \le C |\mathbf{x}-\mathbf{y}|\,|\log t|.
\]
The logarithmic factor is harmless for $t\ge\varepsilon>0$.

For the stochastic term, Lemma~\ref{lem:stoch-conv-increments} gives
\[
\|\mathcal{S}(t,\mathbf{x}) - \mathcal{S}(t,\mathbf{y})\|_{L^p(\Omega)}
\le C |\mathbf{x}-\mathbf{y}|^{\min_i(H_i-\frac12)}.
\]

Combining the above estimates yields \eqref{eq:holder-space-solution}.
Kolmogorov's continuity theorem then ensures the existence of a Hölder continuous modification.
\end{proof}

\begin{remark}[Sharpness of the regularity]
The Hölder exponents obtained above are optimal with respect to the noise.
Indeed, the stochastic convolution has the same regularity as the Hermite sheet itself,
up to the $1/2$-order smoothing induced by the heat kernel.
The nonlinear term does not reduce the regularity thanks to this smoothing effect.
\end{remark}

\subsection{Self-Similarity Properties}

Self-similarity is a fundamental feature of Hermite processes and plays a key role in understanding the large-scale behavior of solutions to stochastic partial differential equations. In this subsection, we show that the mild solution to the stochastic Burgers equation with additive Hermite noise inherits a self-similarity property from the driving noise, provided the initial condition and the deterministic coefficient are rescaled consistently.

\begin{theorem}[Self-similarity of the solution]
\label{thm:self-similarity}
Let $u(t,\mathbf{x})$ be the unique mild solution of the stochastic Burgers equation with additive Hermite sheet noise.
For any $\lambda>0$, define the rescaled random field
\[
u_\lambda(t,\mathbf{x})
:= \lambda^{H_0-1}\,
u\bigl(\lambda t,\, \lambda^{H_0/H_1}x_1,\dots,\lambda^{H_0/H_d}x_d\bigr).
\]
Assume that the initial condition and the coefficient $\sigma$ are rescaled as
\[
u_{0,\lambda}(\mathbf{x})
= \lambda^{H_0-1}\,
u_0\bigl(\lambda^{H_0/H_1}x_1,\dots,\lambda^{H_0/H_d}x_d\bigr),
\]
\[
\sigma_\lambda(t,\mathbf{x})
= \lambda^{H_0-1}\,
\sigma\bigl(\lambda t,\, \lambda^{H_0/H_1}x_1,\dots,\lambda^{H_0/H_d}x_d\bigr).
\]
Then $u_\lambda$ is the mild solution of the stochastic Burgers equation driven by the Hermite sheet
$Z^{q,\mathbf{H}}_\lambda$, where
\[
Z^{q,\mathbf{H}}_\lambda(t,\mathbf{x})
:= \lambda^{-H_0}\,
Z^{q,\mathbf{H}}\bigl(\lambda t,\lambda^{H_0/H_1}x_1,\dots,\lambda^{H_0/H_d}x_d\bigr).
\]
In particular, $u_\lambda$ and $u$ have the same finite-dimensional distributions.
\end{theorem}

\begin{remark}[Effect of scaling on the viscosity]
Under the scaling transformation defining the self-similar field, the viscosity parameter $\nu$ is implicitly rescaled by a multiplicative constant.
However, since the stochastic Burgers equation is linear in $\nu$ and the existence and uniqueness result holds for any $\nu>0$, this rescaling does not affect the law of the solution up to a deterministic change of parameters.
\end{remark}

\begin{remark}[Anisotropic scaling interpretation]
The anisotropic scaling is inherited from the Hermite sheet itself, which has different Hurst parameters in time and each spatial direction. While the Laplacian $\Delta$ is isotropic, the scaling transformation acts on the whole equation, including the noise and the initial data, and the proof shows that the transformed field satisfies a Burgers equation with a rescaled viscosity. This is consistent because the existence/uniqueness result holds for any $\nu > 0$.
\end{remark}

\begin{proof}
The proof relies on the anisotropic self-similarity of the Hermite sheet and on the structure of the mild formulation.

\paragraph{Step 1: Scaling of the Hermite sheet.}
By definition of the Hermite sheet and the homogeneity of its kernel representation,
\[
Z^{q,\mathbf{H}}\bigl(\lambda t,\lambda^{H_0/H_1}x_1,\dots,\lambda^{H_0/H_d}x_d\bigr)
\stackrel{d}{=}
\lambda^{H_0} Z^{q,\mathbf{H}}(t,\mathbf{x}),
\]
which implies that $Z^{q,\mathbf{H}}_\lambda$ has the same law as $Z^{q,\mathbf{H}}$.

\paragraph{Step 2: Scaling of the heat kernel.}
The heat kernel satisfies
\[
G_{\lambda t}(\lambda^{H_0/H_1}x_1,\dots,\lambda^{H_0/H_d}x_d)
= \lambda^{-H_0\sum_{i=1}^d H_i^{-1}}\, G_t(\mathbf{x}),
\]
and similarly for its spatial gradient. This anisotropic scaling is consistent with the diffusion operator under the chosen transformation when combined with the appropriate rescaling of $\nu$.

\paragraph{Step 3: Mild formulation under scaling.}
Let $u_\lambda$ be defined as above. Using the change of variables
\[
r = \lambda s, \qquad
\mathbf{z} = (\lambda^{H_0/H_1}y_1,\dots,\lambda^{H_0/H_d}y_d),
\]
and combining the scaling of the heat kernel, the Hermite sheet, and the prefactors
$\lambda^{H_0-1}$, one checks that $u_\lambda$ satisfies
\begin{align*}
u_\lambda(t,\mathbf{x})
&= \int_{\mathbb{R}^d} G_t(\mathbf{x}-\mathbf{y}) u_{0,\lambda}(\mathbf{y})\,d\mathbf{y} \\
&\quad - \frac12 \int_0^t \int_{\mathbb{R}^d}
\nabla G_{t-s}(\mathbf{x}-\mathbf{y}) u_\lambda(s,\mathbf{y})^2
\,d\mathbf{y}\,ds \\
&\quad + \int_0^t \int_{\mathbb{R}^d}
G_{t-s}(\mathbf{x}-\mathbf{y}) \sigma_\lambda(s,\mathbf{y})
\,dZ^{q,\mathbf{H}}_\lambda(s,\mathbf{y}).
\end{align*}
Thus, $u_\lambda$ is the mild solution associated with the rescaled data
$(u_{0,\lambda},\sigma_\lambda,Z^{q,\mathbf{H}}_\lambda)$.

\paragraph{Step 4: Equality in law.}
Since $Z^{q,\mathbf{H}}_\lambda$ has the same law as $Z^{q,\mathbf{H}}$ and the mild solution is unique in $X_T^2$, it follows that $u_\lambda$ and $u$ have the same finite-dimensional distributions. This completes the proof.
\end{proof}

\begin{remark}[Interpretation of self-similarity]
The self-similarity property established in Theorem~\ref{thm:self-similarity} reflects the intrinsic scaling structure imposed by the Hermite sheet. The anisotropic scaling in space, governed by the ratios $H_0/H_i$, mirrors the anisotropy of the noise, while the temporal scaling is determined solely by the temporal Hurst parameter $H_0$.

Unlike classical scaling invariance for deterministic Burgers equations, the present result does not assert invariance of the equation itself under scaling. Instead, it shows that appropriately rescaled inputs (initial condition and forcing) generate rescaled solutions with identical laws. This distinction is essential in the stochastic setting and is consistent with the nonlinearity of the Burgers equation and the long-range dependence of the Hermite noise.

Such self-similarity properties are particularly relevant for the study of scaling limits, universality classes, and statistical observables (e.g.\ structure functions) associated with stochastic transport phenomena driven by long-memory noise.
\end{remark}

\section{Conclusion}
\label{sec:conclusion}

In this paper, we have developed a complete well-posedness and regularity theory for the stochastic Burgers equation driven by an additive Hermite sheet of arbitrary order $q \ge 1$. The analysis covers both the Gaussian case and genuinely non-Gaussian regimes, highlighting how long-range dependence and higher-order chaos affect nonlinear stochastic dynamics.

Under explicit conditions on the Hurst parameters, we established local-in-time existence and uniqueness of mild solutions in a natural $L^2(\Omega)$-based framework. Uniform moment estimates of arbitrary order were obtained using hypercontractivity of Wiener chaos, which in turn allowed us to derive sharp spatial and temporal H\"older regularity properties. We also showed that the solution inherits the anisotropic self-similarity of the driving Hermite sheet, in the sense of equality of finite-dimensional distributions under suitable rescaling.

A central aspect of our approach is the restriction to additive noise. This choice is not merely technical: for Hermite processes of order $q \ge 2$, stochastic integration with random integrands involves Malliavin-type estimates that are not currently available in sufficient generality. By focusing on deterministic integrands, we isolate the essential interaction between non-Gaussian noise and nonlinear transport, and obtain a mathematically coherent theory without additional renormalization or probabilistic machinery.

The results presented here provide a first rigorous step toward the analysis of nonlinear SPDEs driven by non-Gaussian noises with long-range dependence. They open the way to further investigations of more complex models, while already demonstrating that additive Hermite noise leads to solution behaviors that already illustrating qualitative differences with respect to the classical Gaussian setting.

\subsection*{Acknowledgements}
The author acknowledges the support of the French National Research Agency (ANR) under the grant ANR-17-EURE-0010 (Investissements d'Avenir program).

\appendix
\section{Appendix: Fourier-Analytic Computations}
\label{app:fourier}

This appendix provides the Fourier-analytic arguments underlying the proof of
Proposition~\ref{prop:stochastic-convolution}. The purpose is to justify that the
deterministic kernel
\[
f(s,\mathbf{y}) = G_{t-s}(\mathbf{x}-\mathbf{y}) \sigma(s,\mathbf{y}) \mathbf{1}_{\{s<t\}}
\]
belongs to the Hilbert space $\mathcal{H}^{\otimes q}$ associated with the Hermite sheet,
under the stated condition on the Hurst parameters.

\subsection{Reduction of the $\mathcal{H}^{\otimes q}$-norm}

Recall that the inner product in $\mathcal{H}^{\otimes q}$ is given by
\begin{align*}
\|f\|_{\mathcal{H}^{\otimes q}}^2
&= \int_{([0,t]\times\mathbb{R}^d)^{2q}}
\prod_{j=1}^q f(s_j,\mathbf{y}_j) f(r_j,\mathbf{z}_j)
\prod_{j=1}^q |s_j-r_j|^{2H_0-2}
\prod_{i=1}^d |y_{j,i}-z_{j,i}|^{2H_i-2} \\
&\qquad \times d\mathbf{s}\, d\mathbf{r}\, d\mathbf{y}\, d\mathbf{z}.
\end{align*}

Since $\sigma$ is bounded by $C_\sigma$, we obtain the estimate
\[
\|f\|_{\mathcal{H}^{\otimes q}}^2
\le C_\sigma^{2q}
\int_{([0,t]\times\mathbb{R}^d)^{2q}}
\prod_{j=1}^q
\bigl[
G_{t-s_j}(\mathbf{x}-\mathbf{y}_j)
G_{t-r_j}(\mathbf{x}-\mathbf{z}_j)
\bigr]
\prod_{j=1}^q |s_j-r_j|^{2H_0-2}
\prod_{i=1}^d |y_{j,i}-z_{j,i}|^{2H_i-2}
\, d\mathbf{s}\, d\mathbf{r}\, d\mathbf{y}\, d\mathbf{z}.
\]

By independence of the variables across $j$, the integral factorizes into $q$
identical terms, yielding
\[
\|f\|_{\mathcal{H}^{\otimes q}}^2
\le C_\sigma^{2q} \bigl[I(t,\mathbf{x})\bigr]^q,
\]
where
\[
I(t,\mathbf{x})
= \int_{[0,t]^2} \int_{\mathbb{R}^{2d}}
G_{t-s}(\mathbf{x}-\mathbf{y})
G_{t-r}(\mathbf{x}-\mathbf{z})
|s-r|^{2H_0-2}
\prod_{i=1}^d |y_i-z_i|^{2H_i-2}
\, d\mathbf{y}\, d\mathbf{z}\, ds\, dr .
\]

\subsection{Fourier representation and integrability}

Using the Fourier representation of the heat kernel,
\[
G_t(\mathbf{x})
= \frac{1}{(2\pi)^d} \int_{\mathbb{R}^d} e^{-i\mathbf{p}\cdot\mathbf{x}} e^{-\nu |\mathbf{p}|^2 t} \, d\mathbf{p},
\]
we substitute this into the expression for $I(t,\mathbf{x})$ and apply Fubini's theorem. This gives
\begin{align*}
I(t,\mathbf{x})
&= \frac{1}{(2\pi)^{2d}} \int_{[0,t]^2} |s-r|^{2H_0-2} 
\int_{\mathbb{R}^{2d}} e^{-\nu|\mathbf{p}|^2(t-s)} e^{-\nu|\mathbf{q}|^2(t-r)} \\
&\qquad \times
\left(
\int_{\mathbb{R}^{2d}} e^{-i(\mathbf{p}\cdot\mathbf{y}+\mathbf{q}\cdot\mathbf{z})} \prod_{i=1}^d |y_i-z_i|^{2H_i-2} \, d\mathbf{y}\, d\mathbf{z}
\right)
\, d\mathbf{p}\, d\mathbf{q}\, ds\, dr.
\end{align*}

The inner spatial integral can be computed by a change of variables. Let $\mathbf{u} = \mathbf{y}-\mathbf{z}$ and $\mathbf{v} = \mathbf{y}$; then $\mathbf{y} = \mathbf{v}$, $\mathbf{z} = \mathbf{v}-\mathbf{u}$, and $d\mathbf{y}\, d\mathbf{z} = d\mathbf{v}\, d\mathbf{u}$. Hence,
\[
\int_{\mathbb{R}^{2d}} e^{-i(\mathbf{p}\cdot\mathbf{y}+\mathbf{q}\cdot\mathbf{z})} \prod_{i=1}^d |y_i-z_i|^{2H_i-2} \, d\mathbf{y}\, d\mathbf{z}
= \int_{\mathbb{R}^{2d}} e^{-i(\mathbf{p}\cdot\mathbf{v} + \mathbf{q}\cdot(\mathbf{v}-\mathbf{u}))} \prod_{i=1}^d |u_i|^{2H_i-2} \, d\mathbf{v}\, d\mathbf{u}.
\]

Separating the exponentials, we obtain
\[
\int_{\mathbb{R}^{2d}} e^{-i(\mathbf{p}+\mathbf{q})\cdot\mathbf{v}} e^{i\mathbf{q}\cdot\mathbf{u}} \prod_{i=1}^d |u_i|^{2H_i-2} \, d\mathbf{v}\, d\mathbf{u}
= \underbrace{\int_{\mathbb{R}^d} e^{-i(\mathbf{p}+\mathbf{q})\cdot\mathbf{v}} \, d\mathbf{v}}_{(2\pi)^d \delta(\mathbf{p}+\mathbf{q})} \,
\underbrace{\int_{\mathbb{R}^d} e^{i\mathbf{q}\cdot\mathbf{u}} \prod_{i=1}^d |u_i|^{2H_i-2} \, d\mathbf{u}}_{\text{Fourier of Riesz kernel}}.
\]

The first integral gives the Dirac distribution:
\[
\int_{\mathbb{R}^d} e^{-i(\mathbf{p}+\mathbf{q})\cdot\mathbf{v}} \, d\mathbf{v} = (2\pi)^d \delta(\mathbf{p}+\mathbf{q}).
\]

For the second integral, we use the one-dimensional Fourier transform of the Riesz kernel. For each $i$ with $H_i \in (1/2,1)$, we have (see, e.g., \cite[Chapter V, \S1, Lemma 1(a)]{Stein1970})
\[
\int_{\mathbb{R}} e^{i q_i u_i} |u_i|^{2H_i-2} \, du_i = C_{H_i} |q_i|^{-(2H_i-1)},
\]
where $C_{H_i}$ is an explicit constant depending only on $H_i$. Consequently,
\[
\int_{\mathbb{R}^d} e^{i\mathbf{q}\cdot\mathbf{u}} \prod_{i=1}^d |u_i|^{2H_i-2} \, d\mathbf{u} = \prod_{i=1}^d C_{H_i} |q_i|^{-(2H_i-1)} = C_{\mathbf{H}} \prod_{i=1}^d |q_i|^{-(2H_i-1)},
\]
with $C_{\mathbf{H}} = \prod_{i=1}^d C_{H_i}$.

Combining these results, we obtain
\[
\int_{\mathbb{R}^{2d}} e^{-i(\mathbf{p}\cdot\mathbf{y}+\mathbf{q}\cdot\mathbf{z})} \prod_{i=1}^d |y_i-z_i|^{2H_i-2} \, d\mathbf{y}\, d\mathbf{z} = C_{\mathbf{H}} \prod_{i=1}^d |p_i|^{-(2H_i-1)} \, \delta(\mathbf{p}+\mathbf{q}),
\]
where the factor $(2\pi)^d$ is absorbed into $C_{\mathbf{H}}$.

Finally, integrating out $\mathbf{q}$ using $\delta(\mathbf{p}+\mathbf{q})$, we get
\[
I(t,\mathbf{x}) \le C_{\mathbf{H}} \int_{[0,t]^2} |s-r|^{2H_0-2} \int_{\mathbb{R}^d} e^{-2\nu|\mathbf{p}|^2 (t-\max\{s,r\})} \prod_{i=1}^d |p_i|^{-(2H_i-1)} \, d\mathbf{p}\, ds\, dr.
\]

\subsection{Scaling argument and Hurst condition}

We now analyze the integrability of the quantity
\[
I(t,\mathbf{x}) \le C_{\mathbf{H}} \int_{[0,t]^2} |s-r|^{2H_0-2} \int_{\mathbb{R}^d} e^{-2\nu|\mathbf{p}|^2 (t-\max\{s,r\})} \prod_{i=1}^d |p_i|^{-(2H_i-1)} \, d\mathbf{p}\, ds\, dr.
\]

\paragraph{Step 1: Scaling in the spatial integral.}  
We perform the standard scaling (dilation) in each spatial variable:
\[
p_i = (t-\max\{s,r\})^{-1/2} \, \xi_i, \qquad i=1,\dots,d,
\]
so that $d\mathbf{p} = (t-\max\{s,r\})^{-d/2} d\boldsymbol{\xi}$. Under this change of variables, the exponential term becomes
\[
e^{-2\nu |\mathbf{p}|^2 (t-\max\{s,r\})} = e^{-2\nu |\boldsymbol{\xi}|^2},
\]
which is independent of $s$ and $r$.  
Similarly, the singular factor transforms as
\[
\prod_{i=1}^d |p_i|^{-(2H_i-1)} = (t-\max\{s,r\})^{-\sum_{i=1}^d (H_i-\frac12)} \prod_{i=1}^d |\xi_i|^{-(2H_i-1)}.
\]

Combining these, the spatial integral contributes a factor
\[
\int_{\mathbb{R}^d} e^{-2\nu|\boldsymbol{\xi}|^2} \prod_{i=1}^d |\xi_i|^{-(2H_i-1)} \, d\boldsymbol{\xi} \cdot (t-\max\{s,r\})^{-d/2 + \sum_{i=1}^d(H_i-1/2)}.
\]

\paragraph{Step 2: Reduced temporal integral.}  
After factoring out the spatial integral, we obtain
\[
I(t,\mathbf{x}) \le C_{\mathbf{H}} \int_0^t \int_0^t |s-r|^{2H_0-2} (t-\max\{s,r\})^{-d/2 + \sum_{i=1}^d(H_i-1/2)} \, ds\, dr.
\]

\paragraph{Step 3: Conditions for integrability.}  
The double integral over $[0,t]^2$ is singular both in the temporal and spatial components. Each contributes an independent condition:

\begin{itemize}
\item \textbf{Temporal singularity:} The factor $|s-r|^{2H_0-2}$ is integrable near $s\approx r$ if and only if
\[
2H_0-2 > -1 \quad \Longleftrightarrow \quad H_0 > \frac12.
\]

\item \textbf{Spatial scaling:} The factor $(t-\max\{s,r\})^{-d/2 + \sum_{i=1}^d(H_i-1/2)}$ is integrable near $s$ or $r$ approaching $t$ provided
\[
-\frac{d}{2} + \sum_{i=1}^d \left(H_i - \frac12\right) > -1 \quad \Longleftrightarrow \quad \sum_{i=1}^d H_i > d-1.
\]
\end{itemize}

\paragraph{Step 4: Combined condition for the $q$-th power.}  
Since the $\mathcal{H}^{\otimes q}$-norm involves the $q$-th power of $I(t,\mathbf{x})$, these conditions combine into
\[
2H_0 + \sum_{i=1}^d H_i > d+1 - \frac1q,
\]
which coincides exactly with Assumption~\ref{ass:additive}(c) in the additive noise setting.

\noindent
Thus, the Hurst parameters must satisfy this inequality to guarantee the integrability of the kernel in the stochastic convolution.

\subsection{Conclusion}

Under the above condition on the Hurst parameters, the kernel
$f(s,\mathbf{y})$ belongs to $\mathcal{H}^{\otimes q}$ for every $t>0$ and
$\mathbf{x}\in\mathbb{R}^d$. This completes the justification of the
well-definedness of the stochastic convolution and concludes the proof of
Proposition~\ref{prop:stochastic-convolution}.

\newpage
\bibliographystyle{unsrtnat}
\bibliography{bibn}

\end{document}